%% file: main.tex
\documentclass[11pt]{article}
\usepackage[utf8]{inputenc}
\usepackage[
papersize={5.5in, 8.5in},
left=0.35in,
right=0.35in,
top=0.35in,
bottom=0.5in]{geometry}
\usepackage{moresize}
\usepackage{graphicx}
\usepackage{amsthm}
\usepackage{rotating}
\usepackage{amsfonts,amssymb,amsmath}
\usepackage[usenames]{color}
\usepackage[T1]{fontenc}

\usepackage{longtable}
\usepackage{array}
\usepackage{float}
\usepackage{booktabs}
\usepackage[dvipsnames,usenames,table,xcdraw]{xcolor}

\newcommand{\comment}[1]{}
\newcolumntype{P}[1]{>{\centering\arraybackslash}p{#1}}

\graphicspath{ {images/} }

\title{\bf 
Novel required properties of, and efficient algorithms to seek, perfect cuboids
}
\author{\normalsize 
\textcolor[rgb]{0,0.0,0.0}Aubrey de Grey, Los Gatos, California, USA, aubrey@levf.org \\
Philip Gibbs, Middlesbrough, UK, philegibbs@gmail.com \\
Louie Helm, Bali, Indonesia, louiehelm@protonmail.ch \\
}

\begin{document}

\maketitle

\pagestyle{empty}
\thispagestyle{empty}

\begin{abstract}
We present a novel approach to the age-old question of whether perfect cuboids exist. Our approach consists of two new computer search algorithms, arising from the analysis of "perfect plinths" reported by one of us recently, that are much more efficient than any prior algorithm of which we are aware. Using this approach, we also (i) identify some new two-parameter families of edge cuboids, and (ii) show that large proportions of aspect ratios of faces or internal rectangles cannot be present in a perfect cuboid, a property that was previously reported only for faces with the aspect ratio 4:3.
\end{abstract}

\input{1}

\input{2}

\input{3}

\input{4}
\input{5}
\input{6}

\input{7}

\input{references}
\end{document}

%% file: 1.tex
\section{Definitions}

Following classical literature, we define the following:

- A $\mathbf{Pythagorean}$ $\mathbf{triple}$ (PT) is a triple $(l_1,l_2,h)$ of positive integers such that $l_1^2+l_2^2 = h^2$. A triangle with these values as edge lengths is thus right-angled; $h$ is its $\mathbf{hypotenuse}$ and $l_1$ and $l_2$ are termed the $\mathbf{legs}$. It has been known since Euclid’s time that all PTs can be parameterised as $(2kpq,k(p^2-q^2),k(p^2+q^2))$ for positive integers $(k,p,q)$ with $p,q$ coprime and of opposite parity.

- A rectangle is $\mathbf{Pythagorean}$ iff its edges and diagonal are all of integer length. A cuboid is Pythagorean iff all its faces are Pythagorean. Such a cuboid is also termed an $\mathbf{Euler}$ $\mathbf{brick}$.

- A (hypothetical!) $\mathbf{perfect}$ $\mathbf{cuboid}$ (PC) is an Euler brick whose main diagonal is also of integer length.

- A cuboid is $\mathbf{nearly-perfect}$ (a NPC) iff six of its seven distinct inter-vertex distances (IVDs), i.e. the three edges, the three face diagonals and the main diagonal, are integer length. Thus, an Euler brick is one class of nearly-perfect cuboid, and there are two other classes, termed $\mathbf{face}$ $\mathbf{cuboids}$ and $\mathbf{edge}$ $\mathbf{cuboids}$.

- A hexahedron is a $\mathbf{nearly-cuboid}$ iff all six of its faces are quadrilaterals and its 28 IVDs take between 8 and 13 distinct values. A nearly-cuboid is perfect if all its IVDs are integer length.

- A PT or a cuboid is $\mathbf{primitive}$ iff its edges have a GCD of 1. In the PT case, this is equivalent to the edges being pairwise coprime, or to $k$ being 1 in the parameterisation above.

%% file: 2.tex
\section{Background}

The question whether there exist cuboids whose IVDs are all integers was first explored centuries ago, and interest in it has returned with a vengeance in recent decades; a wide-ranging recent survey of the literature may be found in the first section of [1], so we do not attempt such a survey here.

In recent times, and in the absence of progress towards a proof that no PC exists, it has been natural to seek them computationally; several such attempts are referenced in [1]. However, all have thus far failed.

It is trivial to describe, with three parameters, a cuboid five of whose seven IVDs are integers: one merely defines two PTs of different aspect ratios using Euclid's two-variable parameterisation and scales them so that some leg of one is the same length as some leg of the other. However, no complete parameterisation has been discovered for any of the three types of NPC. Thus, searches for PCs have essentially come down to sifting among examples of a five-integer near miss for cases that satisfy two additional constraints.

Recently, one of us (AdG) reported [2] a preliminary analysis of isosceles rectangular frusta, or "plinths". A cuboid is a special case of a plinth, in which the ceiling is the same size as the base. It turns out that perfect plinths (in which all IVDs are integer) exist; the smallest is shown in Figure 1. In any plinth, the sloping faces are isosceles trapezia; importantly, so are the internal quadrilaterals whose parallel edges are diagonals of the rectangular faces (ACGE and DBFH in Figure 1).

\begin{figure}
\centering
{
\centering
\includegraphics[scale=0.38]{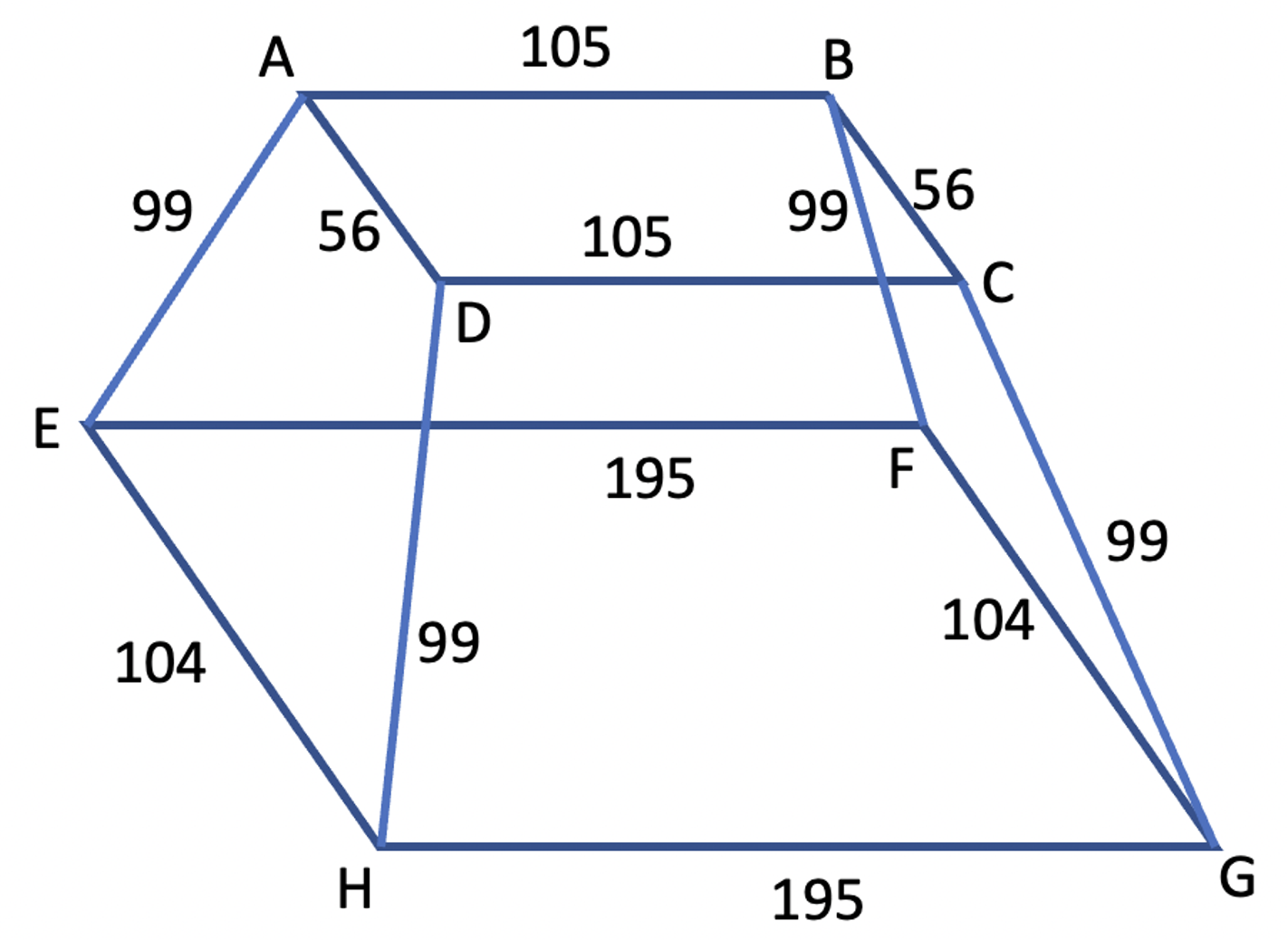} \\
\centering{Figure 1. The smallest perfect plinth. AC=BD=119, EG=FH=221, AH=DE=BG=CF=125, AF=BE=CH=DG=174, AG=BH=CE=DF=190.}
}
\end{figure}

As a result, from any perfect plinth one can generate an infinite family of them, by exploiting the easily-provable fact that the product of the parallel edges of any isosceles trapezium is equal to the square of its diagonal minus the square of its sloping edge. (See Figure 2: $d^2-a^2 = bc$.) Specifically, if for any integers $(i,j)$ we scale the edges of a perfect plinth’s ceiling by $i^2$, those of its base by $j^2$ and its sloping edge by $ij$, the diagonals of all the trapezia will also be scaled by $ij$, so the resulting plinth will also be perfect. Thus, such a family will include a perfect cuboid iff the linear ratio of the ceiling and base of the originating plinth, i.e. $b/c$ of any of the trapezia, is the square of a rational number – or, equivalently, iff $bc$ is a square.

\section{Perfect picture frames}
In the recent report just mentioned [2], a list was given of 15 small perfect plinths. It was overlooked at that time that two of them are degenerate: they have a height of zero, because the difference between the diagonals of the base and the ceiling is exactly twice the "sloping" edge. We term these degenerate perfect plinths "perfect picture frames" (PPFs; see Figure 3). Any family of perfect plinths generated as above always includes a PPF: because degeneracy of the plinth is equivalent to degeneracy of the internal trapezium, the PPF can be derived by scaling (both dimensions of) the originating plinth's base by its internal trapezium’s $c^2$ (i.e. the square of its ceiling’s diagonal), scaling the ceiling by $(d-a)^2$ (the square of the difference between the plinth’s main diagonal and its sloping edge), and scaling the sloping edge by $c(d-a)$, since this also scales the diagonals of the face and internal trapezia by $c(d-a)$. (It can easily be checked that this scaling delivers the degeneracy condition just described.) Thus, any hypothetical perfect cuboid can be derived from a PPF in which the product of the parallel sides of a trapezium (for example, AC times EG in Figure 3) is a square.

\begin{figure}
\centering
{
\centering
\includegraphics[scale=0.38]{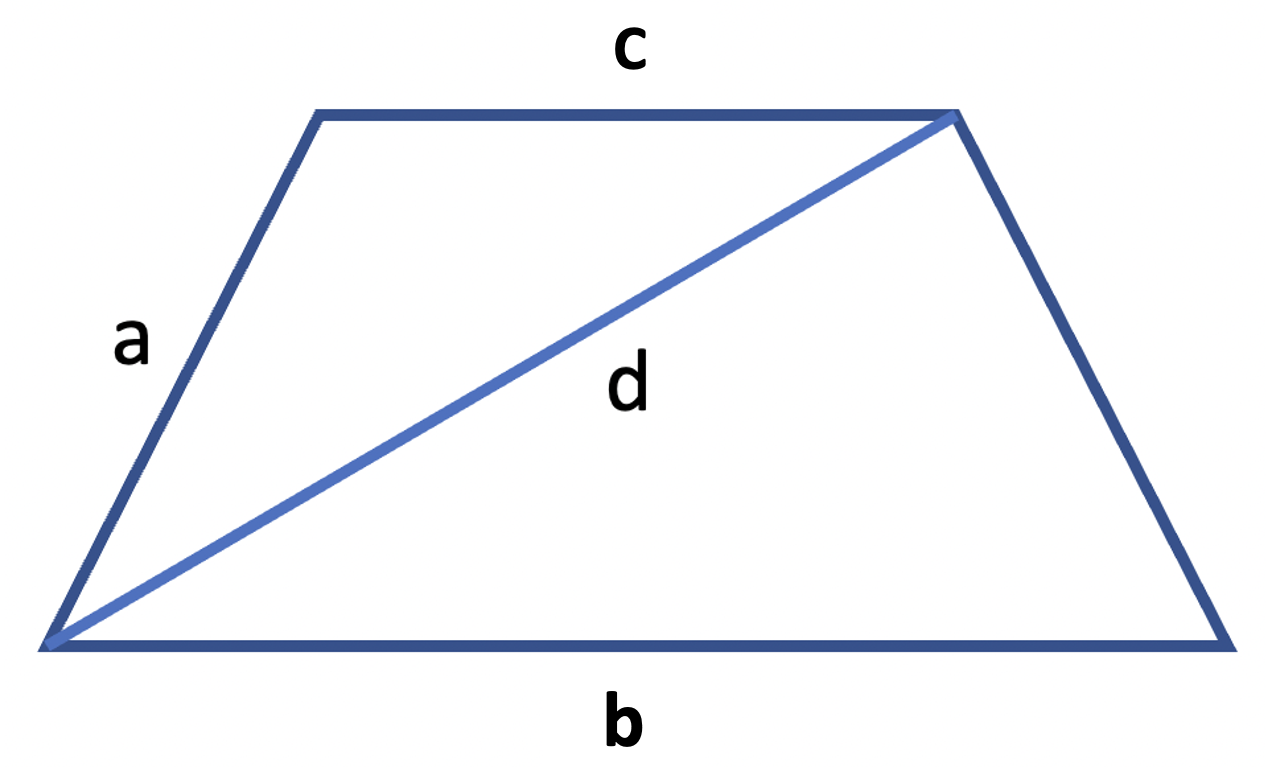} \\
\centering{Figure 2. An isosceles trapezium. In any such trapezium, $d^2-a^2 = bc$.}
}
\end{figure}

PPFs exhibit an additional property that forms the basis for our new search algorithms. As shown in Figure 3, they feature two cases of a set of three PTs (specified in the legend) that share one leg, $k$, and whose other legs $(l,m,n)$ are in geometric progression, i.e. $m=\sqrt{ln}$. (Throughout what follows, we will always assume without loss of generality that $l>n$.) Moreover, the requirement just described for the PPF to correspond to a PC can be expressed in terms of these same numbers: it is just that $l^2-ln$ must be a square. In summary, and noting that $l^2-ln$ is a square iff $ln-n^2$ is, this means that if we could find three numbers $k,l,n$ such that all five quantities $ln, k^2+l^2, k^2+ln, k^2+n^2$ and $l^2-ln$ are squares, we would immediately obtain a PC with edges $\sqrt{ln(k^2+ln)}, l\sqrt{ln-n^2}$ and $k\sqrt{l^2-ln}$, face diagonals $\sqrt{ln(k^2+l^2)}, \sqrt{(k^2+ln)(l^2-ln)}$ and $l\sqrt{k^2+n^2}$ and main diagonal $l\sqrt{k^2+ln}$ – and, conversely, that a proof that no such $(k,l,n)$ exists would constitute a proof that no PC exists.

\begin{figure}
\centering
{
\centering
\includegraphics[scale=0.38]{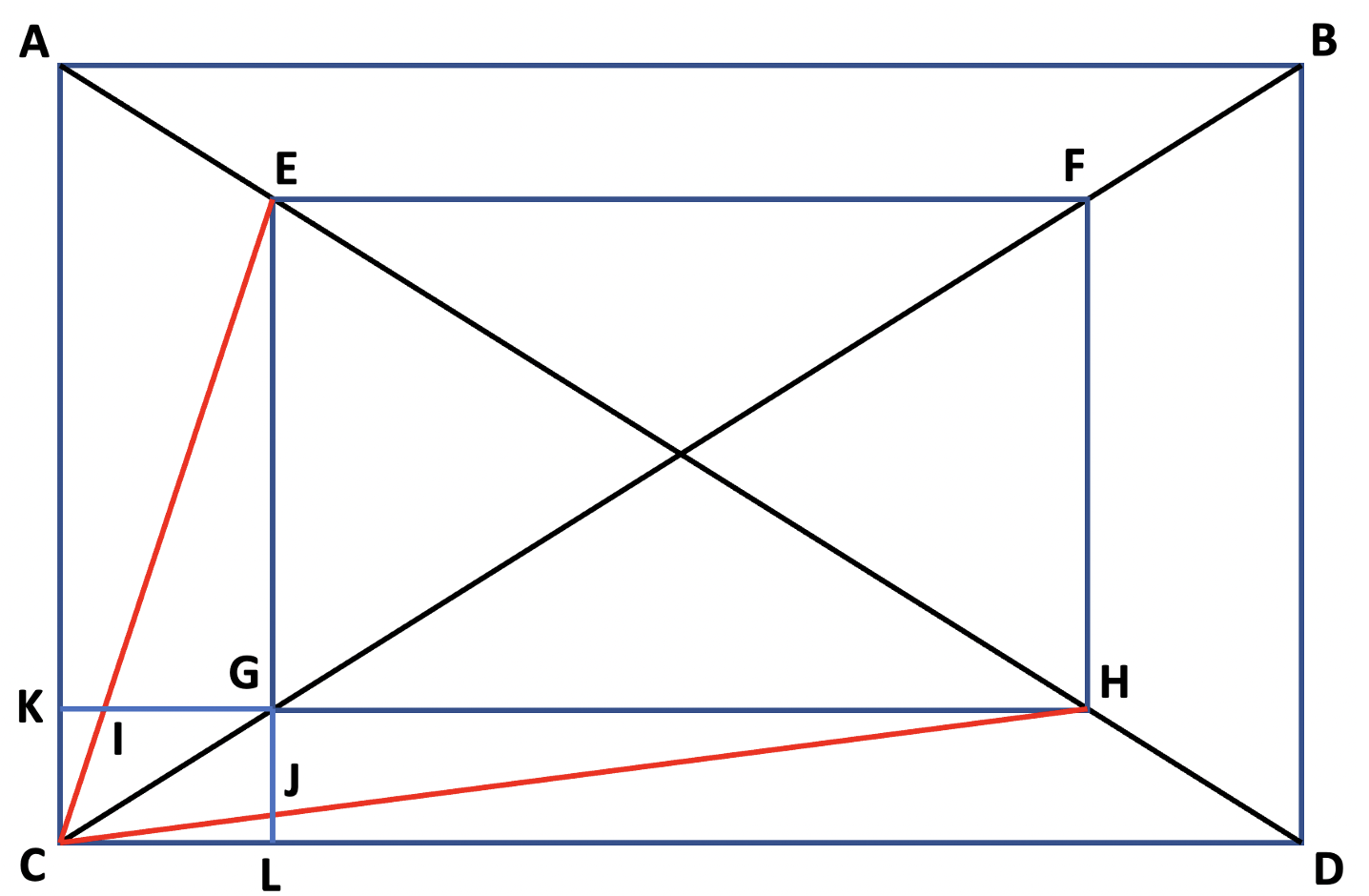} \\
\centering{Figure 3. A “perfect picture frame” (PPF). Perfection holds when all line segments delimited by letters are of integer length. In particular, setting CL = $k$, LE = $l$, LG = $m$ and LJ = $n$, we have $ln=m^2$ and the properties described in the text hold. Similarly for CK, KH, KG and KI.}
}
\end{figure}

\section{Bi-perfect cuboids}
A PPF is one condition short of being the basis of a PC: of the five quantities mentioned above, $ln, k^2+l^2, k^2+ln$ and $k^2+n^2$ are all squares and only $l^2-ln$ is not required to be. We wondered whether similarly interesting families might exist when we instead omit one of the others. Omitting $k^2+l^2$ or $k^2+n^2$ turns out merely to lead to face cuboids. Omitting $k^2+ln$, however, gives a new family. (We also studied the final case, omitting $ln$; it is complex but turned out not to provide further material insight into the PC question, so we summarise our findings in the “Curiosities” section below.) This motivated us to define a generalisation of NPCs that we call “bi-perfect cuboids” (BPCs). In a BPC, the seven IVDs fall into two classes defined by the square-free part (SFP, i.e. an integer divided by its largest square divisor) of their squares: in other words, some of them are squares and all the others have the same SFP. BPCs thus partition into four families:

-	0-BPCs are cuboids all of whose distances are integer, i.e. PCs

-	1-BPCs have exactly one non-integer length, so they are the familiar NPCs, i.e. Euler bricks, face cuboids and edge cuboids

-	2-BPCs have exactly two non-integer lengths and fall into six subclasses (here we treat the case where one of the non-integers is an edge and the other is a face diagonal as two classes, distinguished by whether the two are perpendicular)

-	3-BPCs have exactly three non-integer lengths

There would superficially appear to be several subclasses of 3-BPCs, but in fact it can be shown that all but two of them are impossible, as follows. First consider the seven squares of the IVDs, i.e. $a^2, b^2, c^2, a^2+b^2, a^2+c^2, b^2+c^2$, and $a^2+b^2+c^2$. These will all be positive integers for a 3-BPC; four will be squares and the other three will have the same SFP. Dividing out the greatest common divisor of all seven integers gives seven numbers that fall into a group of four with a SFP of $N$ and a group of three with a SFP of $M$. Call these the $N$-set and the $M$-set respectively. $M$ and $N$ cannot both be 1 (since that would make the cuboid a PC) and must be coprime. Now, suppose $N>1$. If any of the six distinct right-angled triangles in the cuboid has exactly two sides in the $N$-set, then the third side must divide by $N$, since its square is the sum or difference of two numbers that divide by $N$. Thus, if every IVD in the $M$-set is a side of a triangle the other two sides of which are in the $N$-set, all seven distances will divide by $N$, contradicting our earlier removal of all common factors - and this is easily seen to be the case whichever four IVDs comprise the $N$-set. We deduce that $N=1$ and $M>1$. In this setting, when the three IVDs in the $M$-set do not lie in a plane, three of the $N$-set always form a triangle with some pair from the $M$-set and the final one is the sum of all three or the sum of two minus the third, hence similarly divides by $M$, so again we are done. Further, one can also eliminate the $M$-set $(a^2+b^2, a^2+c^2, b^2+c^2)$ by noting that the sums of the quadratic residues of every pair of $a$, $b$ and $c$ modulo any odd prime divisor $p$ of $M$ would all need to be $p$, which can clearly only be the case for at most two of those pairs. This leaves just two cases: where the $M$-set consists of either the sides and diagonal of a face, e.g. $a^2, b^2$ and $a^2+b^2$, or those of an internal rectangle, e.g. $a^2, b^2+c^2$ and $a^2+b^2+c^2$.

Reassuringly, these two are exactly the classes of 3-BPC just described in the discussion of PPFs. Recall that the cuboid arising from a PPF has edges $\sqrt{ln(k^2+ln)}, l\sqrt{ln-n^2}$ and $k\sqrt{l^2-ln}$. Its face diagonals are thus $\sqrt{ln(k^2+l^2)}, l\sqrt{k^2+n^2}$ and $\sqrt{(k^2+ln)(l^2-ln)}$ and its main diagonal is $l\sqrt{k^2+ln}$. Hence, if only $l^2-ln$ is non-integer, we have a 3-BPC in which one face is Pythagorean in aspect ratio and the squares of its sides and diagonal all have the same SFP, while all four remaining lengths are integers. Similarly, if only $k^2+ln$ is non-integer, we have a 3-BPC in which one internal rectangle (whose sides are an edge and a face diagonal of the cuboid, and whose diagonal is the cuboid’s main diagonal) is Pythagorean in aspect ratio and the squares of its sides and diagonal all have the same SFP, while all four remaining lengths are integers. Hereafter we will term these two subclasses “F-BPCs” and “I-BPCs” respectively.

%% file: 3.tex
\section{The first algorithm; congruent number elliptic curves}

\subsection{The utility of square-free parts}
We first note that any pair of integers whose product is a square must have the same SFP: thus, here $SFP(l) = SFP(n)$. Additionally, scaling by any integer $r$ of a PT with legs $a$ and $b$ does not change the SFP of (twice) the PT’s area, i.e. $SFP(ab) = SFP((ra)(rb))$. Accordingly, two PTs with legs $k,l$ and $k,n$, such as CLJ and CLE in Figure 3, must satisfy $SFP(kl) = SFP(kn)$ in order for $\sqrt{ln}$ to be an integer, so they must also satisfy $SFP(k’l’) = SFP(k''n’)$ where $k’$ and $l’$ are coprime, $k’/l’ = k/l$ and correspondingly for $k''$ and $n’$. Conveniently, as noted above, both F-BPCs and I-BPCs have the property that $k^2+l^2$ and $k^2+n^2$ are squares, i.e. the right-angled triangles with legs $(k,l)$ and $(k,n)$ are both Pythagorean.

We can therefore search systematically (and simultaneously) for F-BPCs and I-BPCs by computing and storing $SFP(2pq(p^2-q^2))$ for increasing positive integers $p$ and $q<p$ with $(p,q)$ coprime and of opposite parity. When we encounter a $(p,q)$ whose SFP value we have seen before, we perform the additional checks, with respect to each such earlier $(p’,q’)$, of whether $k^2+ln$ and $l^2-ln$ are squares (where $k:l$ and $k:n$ are $2pq:p^2-q^2$ and $2p’q’:p’^2-q’^2$ in some order), but if the SFP is new we just store it. This algorithm is gratifyingly fast, because it turns out that very few SFP values are common to more than one small $(p,q)$ pair and thus require any additional checks. In fact, among the $(p,q)$ pairs with $p \leq 10^4$, we find that roughly 99.8\% have a SFP that no other has, and this proportion increases as $p$ rises.

\subsection{Congruent number elliptic curves}

The rational points on an elliptic curve form an abelian group, the Mordell-Weil group [3] of the curve. The group factorises into a product of a finite group and an infinite lattice group of the form $\mathbb{Z}^r$ where $r$ is the curve’s rank. The finite part is called the torsion group. If the curve has no rational points other than these torsion points then its rank is zero, otherwise the rank is positive and there are infinitely many rational points on the curve, obtainable from any set of independent generators (of which there are a number equal to the rank) by the well-known tangent and secant methods [3].

Let $n$ be $SFP(pq(p^2-q^2))$ for some Pythagorean parameter pair $(p,q)$, and set $x = np/q$. Now consider the elliptic curve $x^3 - n^2x  = y^2$. Its torsion group consists of the identity (at infinity) and the three trivial  points  $(x,y) = (0,0), (n,0)$ and $(-n,0)$. This equation is known as a congruent number elliptic curve. Any positive integer $n$ for which this curve has rational points other than the torsion points is known as a congruent number. Equivalently, any such $n$ is the area of some Pythagorean triangle (with rational, not necessarily integer, edges).

Therefore, for a fixed SFP there will be either zero non-trivial solutions or infinitely many; we can use freely available algorithms to determine which is the case. When the rank is positive, the two remaining square conditions can be tested for small solutions to determine whether the curve represents a F-BPC or I-BPC (or both, which would mean it represents a PC).

\subsection{Details and results of the first algorithm}

We developed several iterations of software to seek 3-BPCs in this way. We first wrote a Mathematica script on a laptop that located such curves by a brute-force method, of identifying SFP values of $2pq(p^2-q^2)$ for increasing $(p,q)$ and checking for their F-BPC and I-BPC status. This script quickly located the first 78 instances (51 F-BPCs, 27 I-BPCs) with $q<p \leq 8100$. Running a slightly improved version of this Mathematica script on a higher-memory server allowed the search to reach $q<p \leq 81000$ and found 100 additional instances.

Once the Mathematica code was no longer able to efficiently search higher due to memory constraints, we translated its basic algorithm into C++ and located 301 F-BPCs and I-BPCs with $q<p \leq 10^6$.

Next, a more intricate and optimized 700-line C++ version was written which automatically partitioned the search space into 19 separate runs based on the lowest common primes in each SFP (2, 3, 5, 7, etc) while also breaking each run into a 2-pass search. The first pass saved unique signatures of the required parameters into a blocked bloom filter that consumed up to 595 Gb of memory. The second pass queried the bloom filter for paired parameters in parallel on a 72-thread Xeon server with 1.5 Tb of memory. This code was run for several weeks and returned an exhaustive list of 288 F-BPCs and 122 I-BPCs with $q<p \leq 6\cdot10^6$.

Finally, we sought 3-BPCs by analysing the congruent number elliptic curves described above: when the curve has positive rank we derive some rational points on it and ask whether the $(k,l,n)$ arising from two such points equates to a 3-BPC. However, we were initially held back by limitations of the software packages available to us, and then this approach was superseded by our second algorithm (described below). We merely note that the tangent method typically derives PTs with very much larger parameters (as measured by the product of the numerator and denominator) than those of the generator, with the result that 3-BPCs are very unlikely (and we duly found none) for curves whose rank is 1 and thus for which only the tangent method, not the secant method, is initially available. While this report was in review, we became aware of a recent paper [4] in which the authors were able to show that no PC exists corresponding to a congruent number elliptic curve of rank 1.

It goes without saying that we found no example that is both an F-BPC and an I-BPC (or else the title of this report would have been very different!). A notable feature is that, for either class of 3-BPC, the Pythagorean parameters of the rectangle that makes it a 3-BPC (i.e. of $k:\sqrt{ln}$ for F-BPCs or $\sqrt{l^2-ln}:\sqrt{ln}$ for I-BPCs) are always much smaller than those of the originating $k:l$ and $k:n$. Another, to which we briefly return below, is that there are about twice as many F-BPCs as I-BPCs up to a given bound on $p$.

%% file: 4.tex
\section{The “master equation” and new elliptic curves describing BPCs}

We now describe a new method for representing various types of BPCs as elliptic curves, and the implications thereof for the possible shapes of such cuboids.

To find a NPC, the defining equations have previously been reduced to elliptic curves by substituting unknown variables with rational parameterisations of the Pythagorean triples; see e.g. [5]. Here we explore a different approach: multiplying the factors that must be squares to form elliptic curves that split over the rationals, i.e. that take the form $(x+r)(x+s)(x+t) = y^2$. This works because when we find any rational non-torsion points on such an elliptic curve we can use the tangent method to identify rational points at which all three factors $x+r, x+s$ and $x+t$ are squares. Specifically, these elliptic curves have a 4-element torsion group, and a rational solution at $x$ generates three others at $(st - r(x+s+t))/(x+r), (rt - s(x+r+t))/(x+s)$ and $(rs - t(x+r+s))/(x+t)$.

We proceed as follows. Denote the edges of a cuboid as $a, b$ and $c$. Then the following seven factors must be squares in a PC: $a^2, b^2, c^2, a^2+b^2, a^2+c^2, b^2+c^2$ and $a^2+b^2+c^2$. Thus, the product $a^2b^2c^2(a^2+b^2)(a^2+c^2)(b^2+c^2)(a^2+b^2+c^2)$ must also be a square. This ostensibly facile statement turns out to be much more useful than it might initially seem, by providing a basis for the approach just described. We first summarise how it can be used, then we work through the various subtypes in detail.

\subsection{Elliptic curves for 1-BPCs, i.e. NPCs}

In any subtype of NPC, three or four triples among the six integer IVDs form PTs. Thus, in all cases we can specify such a PT (for example we can fix $a^2, b^2$ and $a^2+b^2$ to be squares), and the requirement to be a NPC is then simply that three of the four remaining IVDs must also be integers, i.e. their squares must be squares – which leads to elliptic curves via the property just described.

\subsection{Elliptic curves for 2-BPCs}

When we attempted a similar analysis of the various subtypes of 2-BPCs, we were surprised to find that the resulting curves all turned out to be hyperelliptic rather than elliptic, making it impossible to proceed along the above lines. Perhaps connected to this finding is that we searched computationally for small $(q<p<2^{16})$ examples of each of the six subclasses of 2-BPCs and found none. We were unable, however, to discern a simple explanation for this. 

\subsection{Elliptic curves for 3-BPCs, i.e. F-BPCs and I-BPCs}

It is, however, possible to derive elliptic curves for both of the 3-BPC subtypes. We proceed by treating them as 4-BPCs, i.e. by scaling the cuboid by the square root of the shared SFP. This gives us an equation with four factors, rather than the three in the NPC cases, so not an elliptic curve in the standard Weierstrass form. However, it can readily be transformed into one, as will be described below.

\section{Prohibited, and parametric families of, aspect ratios}

We can now state and prove our main results.

\subsection{Prohibited aspect ratios of faces in a BPC or a PC}

\paragraph{Definitions.} We define the elliptic curves EF1, EF2, EF3 and EF4 as follows:

EF1: $x^3 + (p^2+q^2)^2x^2 + 4p^2q^2(p^2-q^2)^2x = y^2$

EF2: $x^3 + ((p^2-q^2)^4+16p^4q^4)x^2 + 16p^4q^4(p^2-q^2)^4x = y^2$

EF3: $x(x+(p^2-q^2)^2)(x+(p^2-q^2+2pq)^2) = y^2$

EF4: $x(x+4p^2q^2)(x+(p^2-q^2+2pq)^2) = y^2$

\paragraph{Theorem 1.} Let $(p,q)$ be coprime positive integers of opposite parity with $p>q$. Then, if either of the elliptic curves EF1 or EF2, or both of the elliptic curves EF3 and EF4, have rank zero, no PC exists with a face having aspect ratio $(p^2-q^2):2pq$.

\paragraph{Proof.} The proof proceeds by considering various types of BPC, since a PC is a special case of any type of BPC. In particular, considering a PC with integer edges $a,b,c$, and assuming without loss of generality that the largest power of 2 dividing $b$ does not divide $a$, we will examine elliptic curves EF1 through EF4 corresponding to the following cases:

-	$a^2+b^2, a^2+c^2$ and $b^2+c^2$ are squares but $a^2+b^2+c^2$ may be non-square

-	$a^2+b^2, a^2+c^2, b^2+c^2$ and $a^2+b^2+c^2$ may all be non-square, but they all have the same SFP

-	$a^2+b^2, a^2+c^2$ and $a^2+b^2+c^2$ are squares but $b^2+c^2$ may be non-square

-	$a^2+b^2, b^2+c^2$ and $a^2+b^2+c^2$ are squares but $a^2+c^2$ may be non-square

It is easily seen that if the condition on these curves described in the theorem holds, the theorem itself holds.

\paragraph{Remark.} Leech stated [6] that a descent argument (which he seems never to have supplied explicitly) can be used to prove that no PC can have a face with Pythagorean parameters 2 and 1, i.e. with aspect ratio 4:3. Descent arguments are the standard way to show that an elliptic curve has rank zero, i.e. has no rational points other than torsion points. In our curves, torsion points do not correspond to BPCs, so a rank of zero demonstrates that no BPC of the form in question can have a face (resp. internal rectangle) with the specified aspect ratio.

\subsubsection{Case 1: body cuboid (Euler brick), face ratio}

In an Euler brick (that is not a PC) the body diagonal is not rational, i.e. $a^2+b^2+c^2$ is not square. We set $a^2, b^2$ and $a^2+b^2 = d^2$ to be pairwise coprime square integers, or equivalently $a = p^2-q^2, b = 2pq, d = p^2+q^2$. Then we require that (the product of) the remaining three factors $c^2, a^2+c^2$ and $b^2+c^2$ must be square: $c^2(a^2+c^2)(b^2+c^2) = y^2$. Setting $x = c^2$, we then have the elliptic curve $x(x+a^2)(x+b^2) = y^2$, or equivalently (after simplification): \\

$x^3 + (p^2+q^2)^2x^2 + 4p^2q^2(p^2-q^2)^2x = y^2$ \null\hfill (EF1) \\

This curve has an 8-element torsion group with rational points at $x \in {\{0, -a^2, -b^2, ab, -ab\}}$. To form a body cuboid there must be a rational point not in that torsion group – i.e., the rank of the elliptic curve must be positive. The rank (or, at least, bounds on it) can be determined with freely available algorithms. If the rank is found to be zero, there is no Euler brick (and, therefore, no PC) with any face having an aspect ratio $a:b$. Conversely, if the rank is positive, non-torsion points with the $x$-coordinate being a square exist; indeed, any point derived by the tangent method from an arbitrary non-torsion point will have this property.

\subsubsection{Case 2: F-BPC, face ratio}

We again set the face $a^2, b^2, a^2+b^2 = d^2$ to be a Pythagorean triple, or equivalently $a = p^2-q^2, b = 2pq, d = p^2+q^2$. The product of the remaining four factors will then be a square, i.e. $c^2(a^2+c^2)(b^2+c^2)(a^2+b^2+c^2) = w^2$. As it stands, there are four factors here rather than the three in the NPC cases, so this is not an elliptic curve in the usual cubic Weierstrass form. However, it can be transformed into that form with the following steps:

Set $c^2 = 1/z$ giving $(a^2z+1)(b^2z+1)(d^2z+1) = w^2z^4$

Set $w = y/(abdz)^2$ giving $(a^2z+1)(b^2z+1)((d^2z+1) = y^2/(abd)^4$

Set $z = t/(abd)^2$ giving $(t+b^2d^2)(t+a^2d^2)(t+a^2b^2) = y^2$

Set $t = x-a^2b^2$ simplifying to $x(x+a^4)(x+b^4) = y^2$

Combining the above into one rational transformation we can set $x = a^2b^2(a^2+b^2+c^2)/c^2$ and $y = w(a^2b^2(a^2+b^2))^2/c^2$, which gives the elliptic curve $x(x+a^4)(x+b^4) = y^2$, or equivalently: \\

$x^3 + ((p^2-q^2)^4+16p^4q^4)x^2 + 16p^4q^4(p^2-q^2)^4x = y^2$ \null\hfill (EF2) \\

If this curve's rank is positive, the tangent method can be used again to find rational points at which $x, x+a^4$ and $x+b^4$ are all square. Working back from here to the initial four-factor equation we find that $c^2, a^2+c^2, b^2+c^2$ and $a^2+b^2+c^2$ are not necessarily squares but they must all have the same SFP giving the desired F-BPC. Conversely, F-BPCs whose Pythagorean face has aspect ratio $a:b$ cannot exist if this elliptic curve has rank zero, so PCs cannot have such a face either.

\subsubsection{Case 3: face cuboid, face ratio}

Here, one of the squared face diagonals (without loss of generality let it not be $a^2+b^2$) is not a square. Again set $a = p^2-q^2, b = 2pq, d = p^2+q^2$. Then we require that the remaining three factors $c^2, b^2+c^2$ and $a^2+b^2+c^2$, and hence also their product, must be square: $c^2(b^2+c^2)(a^2+b^2+c^2) = y^2$. Setting $x = c^2$, we then have the elliptic curve $x(x+a^2)(x+d^2) = y^2$. The curve in this case therefore has the same structure as for the Euler brick: an 8-element torsion group with rational points at $x \in {\{0, -a^2, -d^2, ad, -ad\}}$. As before, to admit a face cuboid there must be a rational point not in the torsion group, i.e. the rank of the elliptic curve must be positive. If the rank is zero, there is no face cuboid with either of the other faces having an aspect ratio $a:b$.

However, in this setting the equation is not symmetric as between $a$ and $b$, so in order to avoid representing the same aspect ratio by two different $p,q$ pairs (one of them having $p+q$ odd, the other having $pq$ odd) we must distinguish the cases where $a^2+c^2$ is not a square versus where $b^2+c^2$ is not a square. (Recall that we have stipulated, without loss of generality, that the largest power of 2 dividing $b$ does not divide $a$.) This is most straighforwardly achieved by considering two elliptic curves:\\

$x(x+(p^2-q^2)^2)(x+(p^2-q^2+2pq)^2) = y^2$ \null\hfill (EF3)

$x(x+4p^2q^2)(x+(p^2-q^2+2pq)^2) = y^2$ \null\hfill (EF4) \\

If both of these curves have zero rank, no face cuboid (and thus no PC) can have a face with aspect ratio $(p^2-q^2):2pq$. This completes the proof. \qedsymbol

\subsection{First parametric family of edge cuboids}

It will be noted that we have not considered edge cuboids, of which PCs are again a special case. This leads to (and is explained by) our next result.

\paragraph{Theorem 2.} For any coprime positive integers $p$ and $q<p$ of opposite parity satisfying either $p^4-q^4 > 4p^3q$ or $p^4-q^4 < 4pq^3$, there exists an edge cuboid whose integer edges are in the ratio $(p^2-q^2):2pq$.

\paragraph{Proof.} We proceed as in the previous cases; we again write $d^2 = a^2+b^2$ and use the three remaining quantities in a product: $(a^2+c^2)(b^2+c^2)(a^2+b^2+c^2) = y^2$. Setting $x = c^2$ (that is, the square of the irrational edge), this becomes the elliptic curve $(x+a^2)(x+b^2)(x+d^2) = y^2$. However, in contrast to the prior cases, here the curve turns out to possess a rational point (at $x = 0$) that is not a torsion point, regardless of the values of $a, b$ and $d$. The rank of the curve is therefore always positive and an infinite sequence of rational solutions arises via the tangent method, of which some have the $x$-coordinate being the square of a rational. This completes the proof. \qedsymbol

Specifically, for any coprime $(p,q)$ of opposite parity and with $p>q$, we can set the two integer edges of the edge cuboid to be $8p^2q^2(p^4-q^4)$ and $4pq(p^2-q^2)(p^4-q^4)$. Then the three face diagonals are \\

$4pq(p^2-q^2)^2(p^2+q^2)^4$

$|(p^4-4p^2q^2-q^4)(p^4+4p^2q^2-q^4)|$

$|(p^4-2p^3q-2p^2q^2-2pq^3+q^4)(p^4+2p^3q-2p^2q^2+2pq^3+q^4)|$ \\
\\
and the main diagonal is
\\

$(p^4+2p^3q+2p^2q^2-2pq^3+q^4)(p^4-2p^3q+2p^2q^2+2pq^3+q^4)$ \\
\\
while the square of the third edge is
\\

$(p^4-4p^3q-q^4)(p^4+4p^3q-q^4)(p^4-4pq^3-q^4)(p^4+4pq^3-q^4)$ \\

This parameterisation is obtained by applying the tangent method twice to the parameterisation of the curve’s generators, namely $x = (d-a)(d-b) = 2q^2(p-q)^2$ and $y = (d-a)(d-b)(a+b) = x(p^2+2pq-q^2)$; it is the minimal parameterisation, because the parameterised generator itself does not have a squared rational as its $x$-coordinate, while a single application of the tangent method gives a degenerate cuboid with the third edge having length zero. The only further criterion is that $p/q$ must take a value outside the range where $4p^3q \geq p^4-q^4 \geq 4pq^3$, or roughly $1.66 \leq p/q \leq 4.04$, since a value within that range results in the square of the irrational edge being negative. The smallest example, arising from $p=3,q=2$, is the cuboid with edges 7800 and 18720 and main diagonal 24961.

Can an edge cuboid among this family be a PC? We now show that it cannot.

\paragraph{Theorem 3.} There are no coprime positive integers $(p,q<p)$ of opposite parity for which $e = (p^4-4p^3q-q^4)(p^4+4p^3q-q^4)(p^4-4pq^3-q^4)(p^4+4pq^3-q^4)$ is a square.

\paragraph{Proof.} Set $x=(p^2+q^2)^4$ and $y=4p^2q^2(p^2-q^2)^2$. Then $(x-y)^2 - 4xy = e$. Thus, if $e$ is a square, there is an arithmetic progression $(e,(x-y)^2,(x+y)^2$) of squares with step $4xy$. But the step in a 3-element arithmetic progression of squares cannot be a square (by a well-known descent argument), whereas $4xy$ is the square of $4pq(p^2-q^2)(p^2+q^2)^2$. \qedsymbol

\subsection{Prohibited aspect ratios of internal rectangles in a BPC or a PC}

We can analyse the internal rectangles of a hypothetical PC in the same way as above. This leads to our next result.

\paragraph{Definitions.} We define the elliptic curves EI1, EI2, EI3 and EI4 as follows:

EI1: $x(x-(p^2-q^2)^4)(x+8p^2q^2(p^4+q^4)) = y^2$

EI2: $x(x-16p^4q^4)(x+(p^2+q^2)^4-16p^4q^4) = y^2$

EI3: $x^3 - 2(p^4+q^4)x^2 + (p^4-q^4)^2x = y^2$

EI4: $x^3 - ((p^2+q^2)^2+4p^2q^2)x^2 + 4p^2q^2(p^2+q^2)^2x = y^2$

\paragraph{Theorem 4.} Let $(p,q)$ be coprime positive integers of opposite parity with $p>q$. Then if both of the elliptic curves EI1 and EI2, or both of the elliptic curves EI3 and EI4, have rank zero, no PC exists with an internal rectangle having aspect ratio $(p^2-q^2):2pq$.

\paragraph{Remark.} We will not consider Euler bricks here, because by definition no internal rectangle of an Euler brick is Pythagorean (unless the brick is a PC).

\paragraph{Proof.} As with aspect ratios of faces, we proceed by examining the various relevant BPCs.

\subsubsection{Case 1: I-BPC, internal ratio}

For this case we set an internal triangle as a Pythagorean triple, e.g. $a^2, b^2+c^2 = f^2$ and $a^2+b^2+c^2 = a^2+f^2 = g^2$. The product of the remaining four factors will then be a square: $b^2c^2(a^2+b^2)(a^2+c^2) = b^2(b^2+a^2)(f^2-b^2)(g^2-b^2) = w^2$. As with the F-BPC case, there are four factors here rather than the three in the NPC cases, so this is not an elliptic curve in the usual cubic Weierstrass form. However, it can be transformed into one with the following steps:

Set $b^2 = 1/z$ giving $(a^2z+1)(f^2z-1)((a^2+f^2)z-1) = w^2z^4$

Set $w = y/(afgz)^2$ giving $(a^2z+1)(f^2z-1)((a^2+f^2)z-1) = y^2/(afg)^4$

Set $z = t/(afg)^2$ giving $(t+f^2g^2)(t-a^2g^2)(t-a^2f^2) = y^2$

Set $t = x+a^2f^2$ simplifying to $x(x-a^4)(x+f^4+2a^2f^2) = y^2$

Combining the above into one rational transformation we can set  $x = a^2f^2(a^2+f^2-c^2)/c^2$ and $y = w(a^2f^2(a^2+f^2))^2/c^2$; this maps to the elliptic curve $x(x-a^4)(x+f^4+2a^2f^2) = y^2$. Once again the tangent method can be applied when the rank is positive, providing an I-BPC solution. Thus, an I-BPC with internal rectangle aspect ratio $a:f$ cannot exist if this elliptic curve has rank zero. Again we must examine two distinct elliptic curves when we move to the Pythagorean parameter setting, since the defining equation is not symmetric as between $a$ and $f$, and the aspect ratio is only excluded if both curves have zero rank. They are: \\

$x(x-(p^2-q^2)^4)(x+8p^2q^2(p^4+q^4)) = y^2$ \null\hfill (EI1)

$x(x-16p^4q^4)(x+(p^2+q^2)^4-16p^4q^4) = y^2$ \null\hfill (EI2)

\subsubsection{Case 2: face cuboid, internal ratio}

An alternative elliptic curve can be formed for face cuboids, e.g. starting from an internal Pythagorean triangle with sides $c,d,g$. Set $a^2+b^2 = d^2, a^2+b^2+c^2 = d^2+c^2 = g^2$. Then setting $x = a^2$, we have $a^2b^2(b^2+c^2) = x(d^2-x)(g^2-x) = y^2$. Face cuboids cannot, therefore, exist with aspect ratio $d:c$ of either Pythagorean internal rectangle if the corresponding elliptic curve has rank zero. Again we must examine two distinct elliptic curves when working in terms of Pythagorean parameters, since the defining equation is not symmetric as between $c$ and $d$, and the aspect ratio is only excluded if both curves have zero rank. They are, after simplification: \\

$x^3 - 2(p^4+q^4)x^2 + (p^4-q^4)^2x = y^2$ \null\hfill (EI3)

$x^3 - ((p^2+q^2)^2+4p^2q^2)x^2 + 4p^2q^2(p^2+q^2)^2x = y^2$ \null\hfill (EI4) \\

This completes the proof. \qedsymbol

\subsection{Second and third parametric families of edge cuboids}

As in the case of aspect ratios of faces, it will be noted that we have not considered edge cuboids, of which PCs are again a special case. This leads to (and is explained by) our final result.

\paragraph{Theorem 5.} For any coprime positive integers $(p,q)$ of opposite parity with $p>q$ satisfying either $p^4-q^4 < 4p^2q^2$ or $(p^2+q^2)(p-q)^2 > 4p^2q^2$, there exists an edge cuboid with an internal rectangle whose edges are in the ratio $(p^2-q^2):2pq$.

\paragraph{Proof.} We proceed as in the previous cases; we select the internal triangle to be $a^2, b^2+c^2 = e^2$ and $a^2+b^2+c^2 = a^2+e^2 = g^2$. From the product of the remaining factors, proceeding as previously with $x = c^2$ being the square of the irrational edge, we obtain the elliptic curve $b^2(a^2+b^2)(a^2+c^2) = (x+a^2)(e^2-x)(g^2-x) = y^2$. As before, there is a rational point at $x = 0$ that is not a torsion point, and from which the tangent method always generates other rational points whose $x$-coordinates are the squares of rationals, so an edge cuboid exists for any Pythagorean aspect ratio of an internal rectangle. This completes the proof. \qedsymbol

The parameterisation in terms of $a, e$ and $g$ gives generator $x = e^2-ag, y = aeg$. Again we must apply the tangent method twice in order to obtain parameterisations of edge cuboids, because the generator itself does not have a squared rational as its $x$-coordinate and the first application of the tangent method gives a degenerate cuboid with a zero-length edge. We obtain $x = e^2$ - $4e^2a^4g^4/(e^2g^2$ – $a^2e^2 + a^2g^2)^2$, where $x$ is the square of the irrational edge.

When moving to Pythagorean parameters, in this case we have two distinct elliptic curves, since the generator’s $x$-coordinate is not symmetric as between $a$ and $e$; the generators come out as $x = 4p^2q^2-p^4+q^4$ or $(p^2+q^2)(p-q)^2-4p^2q^2$, with $y = 2pq(p^4-q^4)$ in both cases. We then perform two tangent calculations, giving two parameterisations for coprime $(p,q)$ of opposite parity and with $p>q$.

In the case $x = 4p^2q^2-p^4+q^4$, we set an integer edge of the edge cuboid to be
\\

$(p^2-q^2)(p^4+2p^3q+2p^2q^2-2pq^3+q^4)(p^4-2p^3q+2p^2q^2+2pq^3+q^4)$ \\
\\
and the perpendicular face diagonal to be 
\\

$2pq(p^4+2p^3q+2p^2q^2-2pq^3+q^4)(p^4-2p^3q+2p^2q^2+2pq^3+q^4)$ \\
\\
Then the other integer edge is
\\

$4pq(p^4-q^4)^2$ \\
\\
giving the other two face diagonals as 
\\

$(p^2+q^2)(p^4-2p^3q-2p^2q^2-2pq^3+q^4)(p^4+2p^3q-2p^2q^2+2pq^3+q^4)$

$(p^2-q^2)(p^8+8p^6q^2-2p^4q^4+8p^2q^6+q^8)$ \\
\\
and the main diagonal as 
\\

$(p^2+q^2)(p^4+2p^3q+2p^2q^2-2pq^3+q^4)(p^4-2p^3q+2p^2q^2+2pq^3+q^4)$ \\
\\
while the square of the third edge is
\\

$-4 p^2 q^2 (p^4 - 4 p^2 q^2 - q^4) (p^4 + 4 p^2 q^2 - q^4) (3 p^4 + q^4) (p^4 + 3 q^4)$ \\
\\
In this case, values of $p/q$ less than roughly 2.057 give bona fide edge cuboids the square of whose irrational edge is positive; the smallest example, arising from $p=2,q=1$, has integer edges 1443 and 1800 and main diagonal 2405.

In the case $x = (p^2+q^2)(p-q)^2-4p^2q^2$, we set an integer edge of the edge cuboid to be 
\\

$2pq(p^4+2p^3q+2p^2q^2-2pq^3+q^4)(p^4-2p^3q+2p^2q^2+2pq^3+q^4)$ \\
\\
and the perpendicular face diagonal to be 
\\

$(p^2-q^2)(p^4+2p^3q+2p^2q^2-2pq^3+q^4)(p^4-2p^3q+2p^2q^2+2pq^3+q^4)$ \\
\\
Then the other integer edge is
\\

$8p^2q^2(p^2-q^2)(p^2+q^2)^2$ \\
\\
giving the other two face diagonals as
\\

$(p^2+q^2)(p^4-4p^2q^2-q^4)(p^4+4p^2q^2-q^4)$

$2pq(p^8+8p^6q^2-2p^4q^4+8p^2q^6+q^8)$ \\
\\
and the main diagonal is again 
\\

$(p^2+q^2)(p^4+2p^3q+2p^2q^2-2pq^3+q^4)(p^4-2p^3q+2p^2q^2+2pq^3+q^4)$ \\
\\
while the square of the third edge is
\\

$(p^2 - q^2)^2 (p^4 - 2 p^3 q - 2 p^2 q^2 - 2 p q^3 + 
   q^4) (p^4 - 2 p^3 q + 6 p^2 q^2 - 2 p q^3 + q^4) (p^4 + 2 p^3 q - 
   2 p^2 q^2 + 2 p q^3 + q^4) (p^4 + 2 p^3 q + 6 p^2 q^2 + 2 p q^3 + 
   q^4)$ \\
\\
In this case values of $p/q$ greater than roughly 2.89 give bona fide edge cuboids the square of whose irrational edge is positive; the smallest example, arising from $p=4,q=1$, has integer edges 552968 and 554880 and main diagonal 1175057.

%% file: 5.tex
\section{Curiosities of 3-BPCs}

Before discussing our computational analysis of the elliptic curves just described, we note some intriguing features of 3-BPCs that came to our attention during this work.

\subsection{Cousin F-BPCs}

During our investigation of 3-BPCs we noticed that F-BPCs, but not I-BPCs, come in pairs that we term “cousins”. The relationship between a pair is that $2(k^2+ln)/k(l-n)$ of one equals $\sqrt{l/n} - \sqrt{n/l}$ of the other, while they have the same (rational, indeed Pythagorean) value of $k/\sqrt{ln}$. The smallest example is $(k,l,n)$ = (1456,10140,735) and (10032,36100,9801) with $k/\sqrt{ln} = 8/15$. The curves just described for F-BPCs and I-BPCs provide an explanation for this, because they have different torsion groups: the F-BPC curves, but not the I-BPC ones, have a torsion point at $x = a^2b^2$ that generates the cousin transformation. This broadly explains why (as noted above) we found roughly twice as many F-BPCs as I-BPCs up to a given bound on the Pythagorean parameters.

\subsection{Near-PCs with $ln$ not a square}

In the earlier discussion of how to classify BPCs, it was noted that a PC would arise from a triple $k,l,n$ when all five of $ln,k^2+l^2,k^2+ln,k^2+n^2$ and $l^2-ln$ are squares, and that the omission of $l^2-ln$ or $k^2+ln$ gives F-BPCs and I-BPCs respectively. Omitting $k^2+l^2$ or $k^2+n^2$ merely gives face cuboids, but what about omitting $ln$?

We must first take into account that when $ln$ is not a square, we no longer have that $l^2-ln$ is a square iff $ln-n^2$ is. Thus, we must split this case into two: one in which $k^2+l^2,k^2+ln,k^2+n^2$ and $l^2-ln$ are all squares and one in which $k^2+l^2,k^2+ln,k^2+n^2$ and $ln-n^2$ are. Let us call these two cases LM-variants and MN-variants respectively (since we previously used $m$ to denote $\sqrt{ln}$).

It turns out that solutions to both these systems exist, and that they are closely related to F-BPCs and I-BPCs respectively. Let the $k:l$ and $k:n$ ratios of some I-BPC have Pythagorean parameters $(p_1,q_1)$ and $(p_2,q_2)$ and let $h$ be the pairs $(p_1q_2+q_1p_2,|p_1p_2-q_1q_2|)$ with $(p_i,q_i)$ in either order ($p>q$ or $p<q$). Then, there is always a choice of $h$ for which there are two MN-variants, one of whose $(k:l, k:n)$ is shared with the I-BPC and the other is the PT with parameters $h$. The LM-variant case arising from F-BPCs is a little more complex: here we must consider cousin pairs together. If we set $h_1$ and $h_2$ to be defined in the same way for cousin F-BPCs $f_1$ and $f_2$ respectively, we find that four LM-variants arise, with $(k:l, k:n) = (k:l$ of $f_1$, PT with parameters $h_2$) and its counterparts.

We explored these variants further, and identified close connections to other long-standing questions in number theory. However, we did not discover any promising ways forward regarding the PC question.

\subsection{Relationship of F-BPCs to Leech pentacycles}

An intriguing feature of F-BPC cousin pairs is that they appear in the well-known pentacycles described by Leech [6]. For example, one such pentacycle includes the sequence (190/99,15/8,26/7), which corresponds to the F-BPC pair just mentioned: the middle fraction of this trio is the (shared, Pythagorean) aspect ratio $k:\sqrt{ln}$ of the two $(k,l,n)$ tuples, and the fractions flanking it are their $\sqrt{l/n}$. It must be noted that this mapping is not 1-1, i.e. not all pentacycles correspond to F-BPC cousin pairs, because not all pentacycles have a member that is the aspect ratio of a PT (let alone one corresponding to a curve with positive rank).

\subsection{A pair of 3-BPCs with the same SFP}

Although we did not locate a F-BPC + I-BPC pair with a shared SFP, we did locate (with our first algorithm) one F-BPC + F-BPC pair with a shared SFP of 3746496180063, whose elliptic curve has rank 3. The two hits share one $(p,q)$ pair for $(k:l,k:n)$, namely (2678,399); the others are (28798,28779) and (6580798,386019).

This was the only collision of any sort that we encountered, and it may well be the only one that exists, given that the frequency of F-BPCs and I-BPCs dwindles quite rapidly with rising SFP.

%% file: 6.tex
\section{The second algorithm: enumerating 3-BPCs using their elliptic curves}

The above elliptic curves allow us to search for various classes of BPC in an even more efficient manner than with our first algorithm. The search begins, as before, by enumerating increasing $(p,q<p)$ where $(p,q)$ are coprime and of opposite parity. But now, rather than requiring two such $(p,q)$ pairs in order to proceed, we can work with individual ones. Specifically, for the F-BPC case we consider $(p,q)$ to be the Pythagorean parameters of a face, and we examine the elliptic curve that corresponds to the F-BPCs in which that face has integer sides (and diagonal) and the squares of the other four lengths have the same SFP. If that curve’s rank is positive, we obtain generators and compute as many other rational points as we choose using the tangent and secant methods; those for which the $x$-coordinate is a squared rational correspond to F-BPCs. Similarly, for the I-BPC case we assign $(p,q)$ to be the Pythagorean parameters of an internal rectangle and proceed in the same way.

\subsection{Results from the second algorithm}

\subsubsection{F-BPCs and I-BPCs}

This method reproduced the results obtained with our first algorithm, along with 909 additional solutions with larger parameter greater than $6\cdot10^6$, for a total of 1319 hits. We used both Sage [7] and PARI/GP [8], but were handicapped by the lack of robustness of the elliptic curve functions in Sage and the inability of PARI/GP to find generators with large height. We are confident that a more robust implementation of Sage’s algorithm would efficiently enumerate F-BPCs and I-BPCs that are much larger than those we found using our first algorithm. Note that the above does not include arbitrarily large examples derivable using the tangent and secant methods, which are simple to derive from generators.

\subsubsection{Prohibited aspect ratios in PCs}

We also combined the equations in the theorems restricting specific BPC subtypes, to obtain lists of aspect ratios that cannot feature as faces (resp. internal rectangles) of a PC. We did not search to particularly large values of $p$ and $q$ in this regard, but for $p \leq 10^3$ (comprising 202861 admissible $(p,q)$ pairs) we excluded 115637 ratios as faces and 33614 ratios as internal rectangles of any PC. (We are not surprised that the proportion excluded is lower for internal rectangles, since the condition on the curves to test is less stringent.) The simplest aspect ratio not excluded as a face is 21:20 arising from $p=5,q=2$; for internal rectangles it is 15:8 arising from $p=4,q=1$. The proportions of excluded ratios with $p \leq n$ fall with rising $n$ in both cases, but extremely slowly above $n=100$; we speculate, but have not tried to prove, that the excluded ratios retain positive density for unbounded $p$ (maybe even tending to $\frac{1}{2}$ in the face-ratio case).

%% file: 7.tex
\section{Conclusions and directions for future work}

We have developed a new line of attack on the venerable question of whether perfect cuboids exist. It has allowed us to develop an highly efficient algorithm to search for them, which has considerably strengthened the (already strong) heuristic case that none exists. Many avenues are also apparent for taking this approach forward. Among these are:

\subsection{2-BPCs}

We searched for small examples of each of the six subclasses of 2-BPCs and found none. We were unable, however, to find a simple explanation for this, though we suspect it is linked to the fact that our technique for obtaining elliptic curves corresponding to NPCs and 3-BPCs does not work for 2-BPCs. It would be valuable to seek a proof that no 2-BPC exists, because a PC can be considered a special case of a 2-BPC (as it is of any BPC); such a proof may, of course be of the less (but still somewhat!) interesting statement that any 2-BPC must be a PC, but alternatively it may demonstrate the non-existence of PCs.

\subsection{Parameterisations of subsets of edge cuboids}

Saunderson discovered [9], nearly four centuries ago, a family of Euler bricks that is described by just two parameters: for any coprime $(p,q<p)$ of opposite parity, a cuboid with edges 
\\

$8pq(p^4 - q^4)$

$2pq(3p^2 - q^2)(|p^2 - 3q^2|)$ 

$(p^2 - q^2)(|p^2 - 4pq + q^2|)(p^2 + 4pq + q^2)$ \\
\\
has face diagonals 
\\

$2pq(5p^4 - 6p^2q^2 + 5q^4)$

$(p^2 - q^2)(p^4 + 18p^2q^2 + q^4)$ 

$(p^2 + q^2)^3$ \\
\\
It has been shown [10] that no member of this family can be a PC. In the course of our analysis of elliptic curves for BPCs, we found (as described above) some families of edge cuboids with a similar parameterisation, but where the parameters define the aspect ratio of a Pythagorean face or internal rectangle: in other words, there is a member of this family for any such ratio. We showed above that no member of the face-ratio family can be a PC. It would be of interest to prove that members of the two internal-ratio families can also never be PCs – or, conversely, to use those parameterisations as an efficient basis for seeking PCs.

\subsection{Lower bound on the size of a hypothetical PC}

A frustrating limitation of our algorithms is that the failure to find a PC up to a given value of $p$ does not translate into an explicit lower bound on any natural measure of the size of a hypothetical perfect cuboid. This is because the SFP of $2pq(p^2-q^2)$ can be very (arbitrarily?) small relative to $p$ and $q$: for example, it is only 46 when $(p,q)$ = (24336,17689). However, the relative size of the SFP and the parameters of the defining PT (see above) for both classes of 3-BPC leads us to speculate that a subtler lower bound may exist. Our computational analysis using our first algorithm has established the non-existence of PCs with Pythagorean faces or internal rectangles having small parameters, but additionally we found that the shared SFP of the three non-integer IVDs of the 3-BPCs that we found is almost always quite large compared to the Pythagorean parameters of the PT that defines them (namely the PT with aspect ratio $k:\sqrt{ln}$ for F-BPCs or $\sqrt{l^2-ln}:\sqrt{ln}$ for I-BPCs). It would be of great interest if a strong lower bound on some natural measure of the size of a cuboid (such as its volume, the length of its shortest edge, etc) could be established that is a function of those parameters.

\subsection{Near-PCs with $ln$ not a square}

As briefly noted above, we attempted to derive new insights into the PC question by analysing near-PC $(k,l,n)$ triples in which $ln$ is not a square. We did not discover anything over and above our analyses of 3-BPCs, but we may have overlooked a way forward.

\subsection{Relationship of F-BPCs to Leech pentacycles}

We speculate that the intriguing relationship between “cousin” F-BPC pairs and sequences within a Leech pentacycle [6] may have implications that allow further progress in tackling the PC question. This speculation is, however, bereft of any concrete foundation at this point.

\subsection{Implications of the restrictions on aspect ratios}

We have barely scratched the surface of possible consequences of the findings reported here. For example, is there a broad condition $C(p,q)$ such that all elliptic curves corresponding to a particular category of BPC have rank zero when $C(p,q)$ holds? As another example, our computer search up to $p=10^3$ hints that the aspect ratios provably absent from any PC have positive density, but we have not attempted to prove this.